\documentclass[reqno]{amsart}
\usepackage{amssymb}
\usepackage{mathtools}
\usepackage{a4wide,amsmath}
\usepackage{mathrsfs}
\usepackage{amsthm}
\numberwithin{equation}{section}
\numberwithin{figure}{section}
\numberwithin{table}{section}
\usepackage{bbm}
\usepackage{subfig}
\usepackage{enumerate}
\usepackage{needspace}
\usepackage{graphicx}		  
\usepackage{ifpdf}
\ifpdf
\DeclareGraphicsExtensions{.pdf,.eps,.jpg,.png}	
\usepackage[suffix=]{epstopdf}
\fi
\usepackage{xcolor}
\usepackage[utf8]{inputenc}
\usepackage{hyperref}
\hypersetup{hidelinks}

\long\def\MSC#1\EndMSC{\def\arg{#1}\ifx\arg\empty\relax\else
	{\narrower\noindent%
		{2020 Mathematics Subject Classification}: #1\\} \fi}
\long\def\PACS#1\EndPACS{\def\arg{#1}\ifx\arg\empty\relax\else
	{\narrower\noindent%
		{PACS numbers}: #1}\fi}
\long\def\KEY#1\EndKEY{\def\arg{#1}\ifx\arg\empty\relax\else
	{\narrower\noindent%
		Keywords: #1\\}\fi}

\newcommand{\abs}[1]{\lvert#1\rvert} 
\newcommand{\inner}[1]{\langle#1\rangle}

\newcommand{\e}{\mathrm{e}}    
\newcommand{\R}{\mathbb{R}}

\newcommand{\DPhi}{\mathop{\textup{D}\Phi}\nolimits}
\newcommand{\DF}{\mathop{\textup{D}\!F}\nolimits}
\renewcommand{\DH}{\mathop{\textup{D}\!H}\nolimits}
\newcommand{\T}{\textup{T}}

\theoremstyle{plain}
\newtheorem{theorem}{Theorem}[section]

\newtheorem{corollary}[theorem]{Corollary}
\theoremstyle{definition}

\newtheorem{assumption}[theorem]{Assumption}
\theoremstyle{remark}
\newtheorem{remark}[theorem]{Remark}

\begin{document}

\title[Relation between a push-forward and the domain boundary]{Regularity and non-degeneracy of \boldmath{$\Phi^*[I]$} implies regularity of the fixed boundary \boldmath{$\partial\Omega$}}

\author[H.~Garde]{Henrik~Garde}
\address[H.~Garde]{Department of Mathematics, Aarhus University, Aarhus, Denmark.}
\email{garde@math.au.dk}

\author[M.~S.~Vogelius]{Michael~S.~Vogelius}
\address[M.~S.~Vogelius]{Department of Mathematics, Rutgers University, New Brunswick, NJ, USA.}
\email{vogelius@math.rutgers.edu}

\begin{abstract}
	For a diffeomorphism $\Phi$ of a domain $\overline{\Omega}$ onto itself, which is the identity on $\partial\Omega$, we prove that local regularity of the push-forward $\Phi^*[I]$ implies local regularity of $\partial\Omega$, provided a certain non-degeneracy condition is satisfied. To be precise, if $\nu$ is a normal to $\partial\Omega$ at a point $P$, then the condition $(\Phi^*[I](P)-I)\nu \neq 0$ and the assumption that $\Phi^*[I]$ is of class $C^{k+1,\alpha}$ near $P$ imply that $\partial\Omega$ is also of class $C^{k+1,\alpha}$ near $P$. This result naturally complements recent regularity results for non-scattering inhomogeneities.
\end{abstract}

\maketitle



\section{Introduction}

Given a bounded $C^{1,\alpha}$ domain $\Omega\subset \R^n$, $n \ge 2$, and a diffeomorphism $\Phi$ of $\overline{\Omega}$ onto $\overline{\Omega}$, with $\Phi(x)=x$ for $x \in \partial\Omega$ (and $\det \DPhi>0$), the anisotropic matrix-valued coefficient 
\begin{equation*}
	\Phi^*[I]= \frac{\DPhi \DPhi^\T}{\det \DPhi}\circ \Phi^{-1}~,
\end{equation*}
or more generally, the push-forward of a smooth symmetric matrix-valued coefficient $A$,
\begin{equation*}
	\Phi^*[A]= \frac{\DPhi A \DPhi^\T}{\det \DPhi}\circ \Phi^{-1}~,
\end{equation*}
has come to play an important role in the analysis of many electromagnetic imaging problems. To mention a few important facts:
\begin{enumerate}[\rm(i)]
	\item In the context of electrical impedance tomography, media of the form $\Phi^*[A]$ give rise to the same electrical boundary measurements as $A$ (the same Cauchy data)~\cite{KV}, and in two dimensions it is known to be the only media with this property~\cite{ALP}.
	\item In the context of Helmholtz scattering with background conductivity $A$, $\Phi^*[A]$ represent conductivity profiles that are non-scattering 
	for any incident wave.
	\item In the context of the so-called transmission eigenvalue problem, again with ``background conductivity" $A$, $\Phi^*[A]$ represent media for which every value is a transmission eigenvalue.
\end{enumerate}
To understand why these three facts are compatible with generic results in the respective three areas is an interesting challenge. In regards to (i), it is well-known that isotropic conductivities are uniquely determined from the knowledge of all Cauchy data. Therefore the following statement must hold: ``if $A$ and $\Phi^*[A]$ are both isotropic, then $\Phi$ is the identity mapping". Fortunately this is a simple consequence of Liouville's theorem. In regards to (iii), say for background conductivity constant ($\,=I$), very general results, assuming ``well-posedness" of the transmission eigenvalue problem, assert that the spectrum is discrete. Fortunately the push-forward $\Phi^*[I]$ may be easily seen to satisfy the identity
\begin{equation*}
	\inner{\Phi^*[I]\nu,\nu}\inner{\Phi^*[I]\tau,\tau} - \inner{\Phi^*[I]\nu,\tau}^2 = 1 \text{ at all points on } \partial \Omega~, 
\end{equation*}
where $\nu$ is a unit normal vector to $\partial \Omega$ and $\tau$ is any unit tangential vector. This represents a violation of the standard ``covering condition" required for the well-posedness of the transmission eigenvalue problem~\cite{NgNg}. 

This note deals with a compatibility issue in regards to~(ii) and known results to the effect that non-scattering implies a certain regularity of the domain boundary, $\partial\Omega$. The main result in~\cite{CVX} asserts that ``if $A$ is of class $C^{k+1,\alpha}$ (the refractive index is of class $C^{k,\alpha}$) and $\xi^\T(A-I)\nu\ne 0$, where $\xi$ denotes the gradient of a non-scattering incident wave at $P\in \partial \Omega$, then $\partial\Omega$ is of class $C^{k+1,\alpha}$ near $P$". This result has been generalized to domains that are only initially assumed to be Lipschitz in~\cite{KowSaSha}. In order to see the compatibility of these results with~(ii) it therefore becomes necessary to understand why ``if $\Phi^*[I]$ is of class $C^{k+1,\alpha}$ and there exists a vector $\xi$ such that $\xi^\T(\Phi^*[I]-I)\nu \ne  0$, where $\nu$ is the normal to $\partial\Omega$ at a point $P$, then $\partial\Omega$ is necessarily of class $C^{k+1,\alpha}$ near $P$". Expressed differently: ``If $\partial\Omega$ fails to be of class $C^{k+1,\alpha}$ in any neighborhood of a point $P$, then either $(\Phi^*[I](P)-I)\nu=0$ or $\Phi^*[I]$ fails to be of class $C^{k+1,\alpha}$ in any $\overline{\Omega}$-neighborhood of $P$". It is exactly this understanding, which we address in this note. 

It is easy to see that, if $\partial\Omega$ has a proper corner, i.e., is not $C^1$, then any diffeomorphism of $\overline{\Omega}$, $\Phi$, with $\Phi(x)=x$ on $\partial\Omega$, will have $\DPhi = I$ at the corner. In other words, it is easy in this case to understand how $\Phi^*[I]$ violates the non-degeneracy condition of the boundary regularity result. The focus of this note is on domains $\Omega$ that are a priori $C^{1,\alpha}$.

\section{Main results}

Let $\Phi$ be a $C^{1,\alpha}$ diffeomorphism of $\overline{\Omega}$ onto $\overline{\Omega}$ with $\Phi(x) = x$ for all $x\in \partial \Omega$, and suppose $\partial\Omega$ is $C^{1,\alpha}$. Since $\Phi$ is a continuously differentiable invertible mapping on $\overline{\Omega}$, we have that $\det\DPhi \neq 0$ in $\overline{\Omega}$; without loss of generality we shall assume $\det\DPhi >0$. Let $\Phi^*[I](x)$, $x\in\overline{\Omega}$, denote the positive definite symmetric matrix
\begin{equation*}
	\Phi^*[I](x) \coloneqq \frac{\DPhi \DPhi^\T}{\det \DPhi}\circ \Phi^{-1}(x)~.
\end{equation*}
We introduce the following non-degeneracy condition at a point $P\in\partial\Omega$.
\begin{assumption} \label{assump:NonDeg0}
	If $\nu$ is the unit outward normal to $\partial\Omega$ at $P$, assume that
	\begin{equation} \label{NonDeg}
		(\Phi^*[I](P)-I)\nu \neq 0~.
	\end{equation}
\end{assumption}
Stated differently: ``There exists a vector $\xi$ such that $\xi^\T(\Phi^*[I](P)-I)\nu \neq 0$''. Since the above condition is invariant under rotation, we may as well assume $\nu=(0, \dots, 0,1)^\T$. In that case a simple calculation, using the notation $\Phi = (\phi_1,\dots,\phi_n)^\T = (\widetilde{\Phi}^\T,\phi_n)^\T$, gives that
\begin{equation*}
	\DPhi(P) = 
	\begin{bmatrix} 
		\widetilde{I} & \partial_n \widetilde{\Phi} \\
		0      & \partial_n \phi_n 
	\end{bmatrix}\!(P)~,
\end{equation*}
with $\widetilde{I}$ being the  $(n-1)\times(n-1)$ identity matrix, and 
\begin{equation*}
	\det\DPhi(P) = \partial_n \phi_n(P) ~.
\end{equation*}
Therefore 
\begin{equation} \label{eq:Phisimple}
	\Phi^*[I](P) = 
	\begin{bmatrix}
		\bigl\{\frac{\delta_{ij} + \partial_n \phi_i \partial_n\phi_j}{\partial_n\phi_n} \bigr\}_{i,j=1}^{n-1} & \partial_n \widetilde{\Phi} \\[2mm]
		\partial_n \widetilde{\Phi}^\T & \partial_n\phi_n
	\end{bmatrix}\!(P)~.
\end{equation}
The condition \eqref{NonDeg} now translates into
\begin{equation*}
	\partial_n\Phi - (0,\dots,0,1)^\T \neq 0~,
\end{equation*}
at $P$, or in coordinate free notation
\begin{equation} \label{NonDeg3}
	(\DPhi(P)-I) \nu \neq 0~.
\end{equation}
In this paper we prove of the following result.
\begin{theorem} \label{main}
	Suppose Assumption~\ref{assump:NonDeg0} (i.e.,~\eqref{NonDeg3}) is satisfied at a point $P\in\partial\Omega$ and suppose furthermore that $\Phi^*[I]$ is of class $C^{k+1,\alpha}$ near $P$ in $\overline{\Omega}$ for some $k \geq 1$. Then $\partial\Omega$ is of class $C^{k+1,\alpha}$ near~$P$.
\end{theorem}
Since $\Phi^*[I] = \frac{\DPhi \DPhi^\T}{\det \DPhi} \circ \Phi^{-1}$, this theorem immediately implies the following result.
\begin{corollary} \label{main2}
	Suppose Assumption~\ref{assump:NonDeg0} (i.e.,~\eqref{NonDeg3}) is satisfied at a point $P\in \partial \Omega$ and suppose furthermore that $\DPhi\DPhi^\T \circ\, \Phi^{-1}$ is of class $C^{k+1,\alpha}$ near $P$ in $\overline{\Omega}$ for some $k\ge 1$. Then $\partial \Omega$ is of class $C^{k+1,\alpha}$ near~$P$.
\end{corollary}
Due to the formula $\det\Phi^*[I]= (\det\DPhi)^{2-n}\circ \Phi^{-1}$, Theorem~\ref{main} and Corollary~\ref{main2} are actually equivalent for $n \geq 3$.

We note that if Assumption~\ref{assump:NonDeg0} is satisfied at a point $P\in\partial\Omega$ and it is known that $\Phi$ is $C^{k+1,\alpha}$ near $P$ for some $k \geq 1$, then it is fairly easy to show that $\partial\Omega$ is of class $C^{k+1,\alpha}$ near $P$. We also note that, the conditions that $\Phi^*[I]$ or $\DPhi \DPhi^\T \circ\, \Phi^{-1}$ are $C^{k+1,\alpha}$ near $P$ in $\overline{\Omega}$ \emph{do not} guarantee that $\Phi$ is $C^{k+1,\alpha}$ near $P$ in $\overline{\Omega}$. 

Consider the following more restrictive non-degeneracy condition. 
\begin{assumption} \label{assump:NonDeg4}
	Let $\nu$ be the unit outward normal to $\partial\Omega$ at $P$. Assume there exists a vector $\xi$ \emph{tangential to} $\partial\Omega$ at $P$, such that 
	\begin{equation} \label{NonDeg4}
		\xi^\T(\DPhi(P)-I)\nu \neq 0~.
	\end{equation} 
\end{assumption}
We are able to give a very elementary proof of the following analog of Theorem \ref{main}, which only presupposes regularity of $\DPhi \DPhi^\T$ and $\det\DPhi$ on the boundary.
\begin{theorem} \label{main3}
	Suppose Assumption~\ref{assump:NonDeg4} is satisfied at $P\in\partial\Omega$. Select the coordinate system so that $P = 0$ and the $x_n$-axis is in direction of the normal to $\partial\Omega$ at $0$. Assume $\partial\Omega$ near $0$ is given by $x_n = \psi(\tilde{x})$, where $\psi$ locally near $0$ is a $C^{1,\alpha}$ function of $\tilde{x} = (x_1,x_2, \dots , x_{n-1})$. Then we have the following regularity implications:
	\begin{enumerate}[\rm(i)]
		\item If the function $\tilde{x} \mapsto \DPhi \DPhi^\T(\tilde{x},\psi(\tilde{x}))$ is $C^{k,\alpha}$ for some $k\geq 1$ near $0$, then $\partial\Omega$ is of class $C^{k+1,\alpha}$ near $0$.
		\item If the function $\tilde{x} \mapsto \Phi^*[I](\tilde{x},\psi(\tilde{x}))$ is $C^{k,\alpha}$ for some $k\geq 1$ near $0$, then $\partial\Omega$ is of class $C^{k+1,\alpha}$ near $0$.
	\end{enumerate}
\end{theorem}
As we shall see in Section~\ref{sec:proofmain}, Theorem \ref{main3} implies Theorem~\ref{main} (and Corollary~\ref{main2}) if the non-degeneracy condition from Assumption~\ref{assump:NonDeg0} is replaced by Assumption~\ref{assump:NonDeg4}. Hence, the only additional special case is when Assumption~\ref{assump:NonDeg0} is satisfied (c.f.~\eqref{NonDeg3}) but Assumption~\ref{assump:NonDeg4} is not. This is precisely when $\nu$ is an eigenvector of $\DPhi(P)$ for an eigenvalue $\neq 1$. In this special case, we find it necessary to rely on the more technical analysis of boundary regularity for a free boundary problem associated with the anisotropic conductivity $\Phi^*[I]$. 

\section{Proof of Theorem~\ref{main3}}

We first consider the two-dimensional case. Recall that the coordinate system is selected so that $\partial\Omega$ near $P$ (the origin) is given by
$x_2 = \psi(x_1)$ with $\psi'(0) = 0$. With $\Phi = (\phi_1,\phi_2)^\T$, from the fact that $\Phi(x_1,x_2) = (x_1,x_2)^\T$ on $\partial \Omega$, we then obtain that
\begin{align*}
	\partial_1\phi_1[x_1,\psi(x_1)] &= 1 - \partial_2\phi_1[x_1,\psi(x_1)]\psi'(x_1)~, \\
	\partial_1\phi_2[x_1,\psi(x_1)] &= (1- \partial_2\phi_2[x_1,\psi(x_1)])\psi'(x_1)
\end{align*}
near $x_1 = 0$, and thus
\begin{equation*}
	\DPhi[x_1,\psi(x_1)] = \begin{bmatrix} 
		1-\partial_2\phi_1 \psi' & \partial_2\phi_1 \\
		(1- \partial_2\phi_2) \psi'& \partial_2\phi_2
	\end{bmatrix}\![x_1,\psi(x_1)]
\end{equation*}
near $x_1 = 0$. As a consequence, at the boundary point $(x_1,\psi(x_1))$,
\begin{align*}
	\DPhi \DPhi^\T = \begin{bmatrix} 
		(1-\partial_2\phi_1 \psi')^2 + (\partial_2\phi_1)^2 & (1-\partial_2\phi_1\psi')(1-\partial_2\phi_2)\psi' + \partial_2\phi_1 \partial_2\phi_2 \\
		(1-\partial_2\phi_1 \psi')(1-\partial_2\phi_2)\psi' + \partial_2\phi_1\partial_2\phi_2 & (1-\partial_2\phi_2)^2(\psi')^2+(\partial_2\phi_2)^2 
	\end{bmatrix}.
\end{align*}
We now give a detailed proof of Theorem~\ref{main3}(i). Since $x_1 \mapsto \DPhi \DPhi^\T [x_1,\psi(x_1)]$ is of class $C^{k,\alpha}$, it follows that the mapping
\begin{equation} \label{Ck1}
	x_1 \mapsto \begin{bmatrix}
		(1-\partial_2\phi_1\psi')^2 + (\partial_2\phi_1)^2 \\
		(1-\partial_2\phi_1\psi')(1-\partial_2\phi_2)\psi' + \partial_2\phi_1\partial_2\phi_2 \\
		(1-\partial_2\phi_2)^2(\psi')^2 + (\partial_2\phi_2)^2
	\end{bmatrix}\![x_1,\psi(x_1)]
\end{equation}
is of class $C^{k,\alpha}$ near $x_1 = 0$. Now $\det \DPhi(0,0)\neq 0$ and Assumption~\ref{assump:NonDeg4} translate into $\partial_2\phi_1(0,0) \neq 0$ and $\partial_2\phi_2(0,0)\neq 0$. If we introduce the mapping $F\colon \mathbb{R}^3 \rightarrow \mathbb{R}^3$ by
\begin{equation*}
	F(a,b,c) \coloneqq \begin{bmatrix}
		(1-ba)^2 + b^2 \\
		(1-ba)(1-c)a+bc \\
		(1-c)^2a^2+c^2
	\end{bmatrix},
\end{equation*}
then \eqref{Ck1} asserts that 
\begin{equation} \label{Ck2}
	x_1 \mapsto F(\psi'(x_1),\partial_2\phi_1[x_1,\psi(x_1)],\partial_2\phi_2[x_1,\psi(x_1)]) 
\end{equation}
is of class $C^k$ near $x_1=0$. We calculate
\begin{equation*}
	\DF(a,b,c) = \begin{bmatrix}
		-2b(1-ba)    & -2a(1-ba)+2b & 0 \\
		(1-2ba)(1-c) & -(1-c)a^2+c  & -(1-ba)a+b \\
		2a(1-c)^2    & 0            & -2(1-c)a^2 +2c
	\end{bmatrix}.
\end{equation*}
At $(0,b,c)$ this gives
\begin{equation*}
	\DF(0,b,c) = \begin{bmatrix}
		-2b   & 2b & 0 \\
		1-c & c  & b \\
		0     & 0  & 2c
	\end{bmatrix},
\end{equation*}
which is invertible exactly when $bc\neq 0$. By inversion of $F$ it follows from  \eqref{Ck2} that
\begin{equation}
	x_1 \mapsto (\psi'(x_1), \partial_2\phi_1[x_1,\psi(x_1)], \partial_2\phi_2[x_1,\psi(x_1)])^\T
\end{equation}
is of class $C^k$ near $x_1=0$. As a consequence, $\psi$ and therefore $\partial\Omega$ is of class $C^{k+1,\alpha}$ near $0$. This completes the proof of the first regularity implication in the two dimensional case.

In the $n$-dimensional case, suppose the boundary near $P$ (the origin) is given by $x_n = \psi(\tilde{x})$, with $\widetilde{\nabla}\psi(0) = 0$ and $\tilde{x} = (x_1,x_2, \dots,x_{n-1})$. Then for points on $\partial \Omega$, we calculate
\begin{equation*}
	\DPhi = \begin{bmatrix}
		\widetilde{I} - \partial_n\widetilde{\Phi}\widetilde{\nabla}\psi^\T & \partial_n \widetilde{\Phi} \\[1mm]
		(1-\partial_n\phi_n)\widetilde{\nabla}\psi^\T & \partial_n\phi_n
	\end{bmatrix},
\end{equation*}
and thereby
\begin{align*}
	\DPhi\DPhi^\T &= \\
	&\hspace{-1cm}\begin{bmatrix}
		(\widetilde{I} - \partial_n\widetilde{\Phi}\widetilde{\nabla}\psi^\T)(\widetilde{I} - \partial_n\widetilde{\Phi}\widetilde{\nabla} \psi^\T )^\T + \partial_n\widetilde{\Phi}\partial_n\widetilde{\Phi}^\T & (1-\partial_n\phi_n)(\widetilde{\nabla}\psi - \abs{\widetilde{\nabla}\psi}^2 \partial_n\widetilde{\Phi}) + \partial_n\phi_n\partial_n\widetilde{\Phi} \\
		(1-\partial_n\phi_n)(\widetilde{\nabla}\psi^\T - \abs{\widetilde{\nabla}\psi}^2 \partial_n\widetilde{\Phi}^\T) + \partial_n\phi_n\partial_n\widetilde{\Phi}^\T & (1-\partial_n\phi_n)^2\abs{\widetilde{\nabla}\psi}^2 + (\partial_n\phi_n)^2 
	\end{bmatrix}~.
\end{align*}
We now consider the map $F\colon \R^{n-1} \times \R^{n-1} \times \R \to \R^{(n-1)^2+n}$ given by
\begin{equation*}
	F(\tilde{\eta}, \tilde{\xi}, \xi_n) \coloneqq
	\begin{bmatrix}
		(\widetilde{I} - \tilde{\xi}\tilde{\eta}^\T)(\widetilde{I} - \tilde{\eta}\tilde{\xi}^\T) + \tilde{\xi}\tilde{\xi}^\T \\[1mm]
		(1-\xi_n)(\tilde{\eta}^\T - \abs{\tilde{\eta}}^2\tilde{\xi}^\T) + \xi_n\tilde{\xi}^\T \\[1mm]
		(1-\xi_n)^2\abs{\tilde{\eta}}^2 + \xi_n^2
	\end{bmatrix}.
\end{equation*}
For $n=2$ this is the same map $F$ defined earlier. Note that for $n\geq 3$, $(n-1)^2+n$ is strictly greater than $2n-1$. The linearization of $F$ at $(0,\tilde{\xi},\xi_n)$, in direction $(\tilde{k},\tilde{h}, h_n)$, is given by
\begin{equation*}
	\DF(\tilde{k},\tilde{h},h_n) = \begin{bmatrix}
		(\tilde{h}-\tilde{k})\tilde{\xi}^\T + \tilde{\xi}(\tilde{h}-\tilde{k})^\T \\[1mm]
		\xi_n(\tilde{h}-\tilde{k})^\T + \tilde{k}^\T + h_n\tilde{\xi}^\T \\[1mm]
		2\xi_n h_n
	\end{bmatrix}.
\end{equation*}
We note that if $\tilde{\xi}\neq 0$ and $\xi_n\neq 0$ then the linear map $\DF$ is \emph{injective}. To see this, suppose $\DF(\tilde{k},\tilde{h},h_n) = 0$, then it follows immediately from right-multiplication of the first component of $\DF$ by $\tilde{\xi}/\abs{\tilde{\xi}}^2$ and left-multiplication by $\tilde{\xi}^\T/\abs{\tilde{\xi}}$ that $(\tilde{h}-\tilde{k}) \cdot \tilde{\xi} = 0$. Using this fact, right-multiplication of the first component of $\DF$ by $\tilde{\xi}$ immediately gives that $\tilde{h} -\tilde{k} =0$. From the third component of $\DF = 0$ we conclude that $h_n = 0$ and the second component now leads to $\tilde{k} = 0$, in other words $\tilde{h} = \tilde{k} = 0$ and $h_n = 0$. 

The injectivity of $\DF$ is sufficient to guarantee the local existence of a $C^\infty$ map $F^{-1}$  (on the image of $F$) with $F^{-1}\circ F = \mathop{\textup{Id}}$. Now consider $\tilde{\eta} = \widetilde{\nabla}\psi$, $\tilde{\xi} = \partial_n\widetilde{\Phi}$, and $\xi_n = \partial_n\phi_n$. Just as in the two-dimensional case, this shows that $\widetilde{\nabla}\psi$ is $C^{k,\alpha}$ if $\DPhi\DPhi^\T$ is in $C^{k,\alpha}$ near $0$. This verifies Theorem~\ref{main3}(i) for a general~$n$.

Theorem~\ref{main3}(ii) can be proven along the same lines. As already pointed out earlier, it actually follows from Theorem~\ref{main3}(i) if $n\geq 3$.

\section{Proof of Theorem~\ref{main}} \label{sec:proofmain}

To see that Theorem~\ref{main3} implies Theorem~\ref{main} (and Corollary~\ref{main2}) with the non-degeneracy condition from Assumption~\ref{assump:NonDeg0} replaced by Assumption~\ref{assump:NonDeg4}, we proceed as follows: Suppose that $\DPhi \DPhi^\T\circ\, \Phi^{-1}$ or $\Phi^*[I]$ is of class $C^{m+1,\alpha}$ in $\overline{\Omega}$ near $P$ for $m=1$, then $\DPhi \DPhi^\T\circ\, \Phi^{-1}(\tilde{x}, \psi(\tilde{x}))$ or $\Phi^*[I](\tilde{x}, \psi(\tilde{x}))$ is of class $C^{m,\alpha}$ near $P$ (the origin in the setting of Theorem~\ref{main3}). Here we used that $\psi$ is of class $C^{m,\alpha}$ for $m=1$. From Theorem~\ref{main3} it now follows that $\psi$ and therefore $\partial\Omega$ are of class $C^{m+1,\alpha}$ for $m=1$ near $P$. We may use the previous argument recursively for $m$ up to $k$, leading to a proof of Theorem~\ref{main} (or Corollary~\ref{main2}).

Up to a rotation aligning the $x_n$-axis with the normal vector, the Jacobian $\DPhi(P)$ is given by
\begin{equation*}
	\DPhi(P) = \begin{bmatrix}
		\widetilde{I} & \partial_n\widetilde{\Phi} \\
		0 & \partial_n\phi_n
	\end{bmatrix}\!(P)~. 
\end{equation*}
Assumption~\ref{assump:NonDeg0} (c.f.~\eqref{NonDeg3}) therefore amounts to
\begin{equation*}
	\text{either} \quad \partial_n\widetilde{\Phi}(P)\neq 0 \quad \text{or} \quad \partial_n\phi_n(P) \neq 1~. 
\end{equation*}
The only case of Theorem~\ref{main}, for which we have not yet established that $\partial\Omega$ is $C^{k+1,\alpha}$ near $P$, is when Assumption~\ref{assump:NonDeg0} is satisfied but Assumption~\ref{assump:NonDeg4} is not, i.e.,  when
\begin{equation*}
	\DPhi(P) = \begin{bmatrix}
		\widetilde{I} & 0 \\
		0 & \partial_n\phi_n
	\end{bmatrix}\!(P)
	\quad \text{and} \quad \partial_n\phi_n(P)\neq 1~.
\end{equation*}
We proceed to examine this case, thus completing the proof of Theorem~\ref{main}. This part of our proof is based on the use of the Hodograph transform (as in \cite{KindNir}) and it may thus be seen as a simplified version of the proof of Theorem 2.1 in \cite{CVX}.

For simplicity of notation we take $n=2$, with standard basis vectors $\e_1$ and $\e_2$ for $\R^2$. The higher dimensional case follows in exactly the same way. Set $y = \Phi(x)$ and define 
\begin{equation*}
	v(y) \coloneqq x_2-y_2 = (\Phi^{-1}(y))_2 - y_2~.
\end{equation*}
This function satisfies
\begin{align}
\label{veq}
	\nabla_y \cdot (\Phi^*[I]\nabla_y v) &= -\nabla_y \cdot (\Phi^*[I] \e_2) \text{ in } \Omega~,\\
	\nu\cdot(\Phi^*[I]\nabla_y v) &= \nu\cdot(I-\Phi^*[I])\e_2 \text{ on } \partial\Omega~, \text{ and}\\
	v &= 0 \text{ on } \partial\Omega~.
\end{align}
Since $\Phi^*[I]$ is in $C^{k+1,\alpha}$ and $\partial\Omega$ is $C^{1,\alpha}$, then elliptic regularity gives that $v$ lies in $C^{1,\alpha}$. We have that $\nu$ at $P$ equals $\pm \e_2$. Due Assumption~\ref{assump:NonDeg0} and the fact that $\partial_2\phi_1(P) = 0$, we also have that $\e_2\cdot(I-\Phi^*[I])\e_2 =1-\partial_2\phi_2(P) \neq 0$. From the conormal boundary condition for $v$ it now follows that  $\partial_2 v(P) = \frac{1-\partial_2\phi_2(P)}{\partial_2\phi_2(P)} \neq 0$. Let us suppose $\partial_2 v(P) > 0$. Without loss of generality, we may assume $P = (0,0)$. We introduce the mapping
\begin{equation*}
	H \colon (y_1,y_2) \mapsto (z_1,z_2)=(y_1,v(y_1,y_2))~,
\end{equation*}
which maps an $\Omega$-neighborhood of the origin bijectively onto the half-disk
\begin{equation*}
	\{\, (z_1,z_2): z_2 > 0 \text{ and } z_1^2 + z_2^2 \leq r^2 \,\}
\end{equation*}
for $r>0$ sufficiently small. The inverse of this map has the form 
\begin{equation*}
	H^{-1} \colon (z_1,z_2) \mapsto (z_1,w(z_1,z_2))~,
\end{equation*}
where $w$ is also $C^{1,\alpha}$. We now derive a non-linear boundary value problem satisfied by the function $w$ --- see \eqref{weq2}--\eqref{wbc2}. It is from this elliptic PDE problem that we will infer regularity for $w$ and thereby regularity for $\partial\Omega$ near $P$ (through the mapping $H^{-1}$). 

First we calculate
\begin{equation*}
	\DH = \begin{bmatrix}
		1 & 0 \\ 
		\partial_1 v & \partial_2 v
	\end{bmatrix}
\end{equation*}
and so
\begin{equation*}
	\textup{D}(H^{-1})\circ H = (\DH)^{-1} = \begin{bmatrix}
		1 & 0 \\ 
		-\frac{\partial_1 v}{\partial_2 v} & \frac{1}{\partial_2 v}
	\end{bmatrix}.
\end{equation*}
From this we conclude that
\begin{equation*}
	\frac{\partial v}{\partial y_1}=-\frac{\partial w}{\partial z_1} \Big/ \frac{\partial w}{\partial z_2}~, \quad \frac{\partial v}{\partial y_2} = 1 \Big/\frac{\partial w}{\partial z_2}~,
\end{equation*}
and that
\begin{equation*}
	\frac{\partial}{\partial y_1}= \frac{\partial}{\partial z_1} - \Bigl(\frac{\partial w}{\partial z_1}\Big/\frac{\partial w}{\partial z_2}\Bigr) \frac{\partial}{\partial z_2}~, \quad \text{and} \quad  \frac{\partial}{\partial y_2} = \Bigl(1\Big/\frac{\partial w}{\partial z_2}\Bigr) \frac{\partial}{\partial z_2}~.
\end{equation*}
With the notation 
\begin{equation*}
	\widetilde{\nabla}_z w = \frac{1}{\partial w/\partial z_2}
	\begin{bmatrix}
		-\frac{\partial w}{\partial z_1} \\
		1
	\end{bmatrix},
\end{equation*}
these may be rewritten as
\begin{equation*}
	\nabla_y v = \widetilde{\nabla}_z w \quad \text{and} \quad \nabla_y = \e_1
	\frac{\partial}{\partial z_1} + \widetilde{\nabla}_z w \frac{\partial}{\partial z_2}~.
\end{equation*}
Let $\Phi^*[I]_H$ denote the matrix-valued function 
\begin{equation*} 
	\Phi^*[I]_H = \Phi^*[I] \circ H^{-1}~.
\end{equation*}
We thus calculate
\begin{align} \label{rhs}
	\nabla_y \cdot (\Phi^*[I]\nabla_y v) &= \frac{\partial}{\partial z_1}(\e_1\cdot \Phi^*[I]_H \widetilde{\nabla}_z w) + \widetilde{\nabla}_z w \cdot \frac{\partial}{\partial z_2}(\Phi^*[I]_H\widetilde{\nabla}_z w) \notag\\
	&= \frac{\partial}{\partial z_1}(\e_1\cdot \Phi^*[I]_H\widetilde{\nabla}_z w) + \frac{1}{2}\frac{\partial}{\partial z_2}(\widetilde{\nabla}_z w \cdot \Phi^*[I]_H\widetilde{\nabla}_z w) \notag\\
	&\phantom{={}} + \frac{1}{2} \widetilde{\nabla}_z w \cdot \Bigl( \frac{\partial}{\partial z_2} \Phi^*[I]_H\Bigr) \widetilde{\nabla}_z w~.
\end{align}
For the last equality we used the simple fact that 
\begin{equation*}
	\frac{1}{2} \frac{\partial}{\partial z_2}( \widetilde{\nabla}_z w \cdot \Phi^*[I]_H \widetilde{\nabla}_z w ) + \frac{1}{2}\widetilde{\nabla}_z w \cdot \Bigl(\frac{\partial}{\partial z_2} \Phi^*[I]\Bigr)\widetilde{\nabla}_z w = \widetilde{\nabla}_z w \cdot \frac{\partial}{\partial z_2} (\Phi^*[I]_H \widetilde{\nabla}_z w)~.
\end{equation*}
We similarly calculate
\begin{equation*}
	\nabla_y \cdot (\Phi^*[I] \e_2) = \frac{\partial}{\partial z_1}( \e_1\cdot \Phi^*[I]_H \e_2) + \widetilde{\nabla}_z w \cdot\Bigl( \frac{\partial}{\partial z_2}\Phi^*[I]_H \e_2\Bigr)~,
\end{equation*}
and after combination with \eqref{veq} and \eqref{rhs} this gives
\begin{align} \label{weq}
	\frac{\partial}{\partial z_1}\bigl(\e_1\cdot \Phi^*[I]_H(\widetilde{\nabla}_z w + \e_2)\bigr) &+ \frac{1}{2}\frac{\partial}{\partial z_2}(\widetilde{\nabla}_z w \cdot \Phi^*[I]_H\widetilde{\nabla}_z w) \notag\\
	&+ \frac{1}{2} \widetilde{\nabla}_z w \cdot \Bigl( \frac{\partial}{\partial z_2} \Phi^*[I]_H\Bigr)(\widetilde{\nabla}_z w + 2\e_2)= 0~,
\end{align}
for $z$ in a neighborhood of $(0,0)$ in the upper positive half-plane.

Since $\partial\Omega$ is a level set for the original function $v$, we have that $\nu = \pm \nabla_y v/\abs{\nabla_y v}$, so the conormal boundary condition for $v$ leads to the following boundary condition for $w$:
\begin{equation} \label{wbc}
	\widetilde{\nabla}_z w \cdot \Phi^*[I]_H \widetilde{\nabla}_z w + \widetilde{\nabla}_z w \cdot (\Phi^*[I]_H - I)\e_2 = 0~,
\end{equation}
at $z_2 = 0$. The boundary value problem \eqref{weq}--\eqref{wbc} is of the form
\begin{equation} \label{weq2}
	\frac{\partial}{\partial z_1} a_1(z,w,\nabla w) + \frac{\partial}{\partial z_2} a_2(z,w,\nabla w) + a_0(z,w,\nabla w) = 0 \text{ for } z_2 > 0
\end{equation}
with
\begin{equation} \label{wbc2}
	b(z,w,\nabla w) = 0 \text{ for } z_2 = 0~.
\end{equation}
Here
\begin{align*}
	a_1(z,w,\nabla w) &= \e_1\cdot\Phi^*[I]_H(\widetilde{\nabla}_z w + \e_2)~, \\
	a_2(z,w,\nabla w) &= \frac{1}{2} \widetilde{\nabla}_z w \cdot \Phi^*[I]_H\widetilde{\nabla}_z w~, \\
	a_0(z,w,\nabla w) &=\frac{1}{2} \widetilde{\nabla}_z w \cdot \Bigl(\frac{\partial}{\partial z_2} \Phi^*[I]_H\Bigr)(\widetilde{\nabla}_z w + 2\e_2)~, \text{ and } \\
	b(z,w,\nabla w) &= \widetilde{\nabla}_z w \cdot \Phi^*[I]_H \widetilde{\nabla}_z w + \widetilde{\nabla}_z w \cdot (\Phi^*[I]_H-I)\e_2~.
\end{align*}
In this context we note that $\widetilde{\nabla} w$ is a function of $\nabla w$, that $\Phi^*[I]_H$ is a function of $z_1$ and $w$, and that $\frac{\partial}{\partial z_2} \Phi^*[I]_H$ is a function of $z_1$, $w$, and $\nabla w$. More precisely 
\begin{equation*}
	\frac{\partial}{\partial z_2} \Phi^*[I]_H = \Bigl(\frac{\partial}{\partial y_2} \Phi^*[I]\Bigr)_H \frac{\partial w}{\partial z_2}~.
\end{equation*}
We proceed to calculate the principal (highest order terms) of the linearization of the operators on the left-hand side of \eqref{weq2}--\eqref{wbc2} around a given function $w_0$. The operator $\widetilde{\nabla}$ has the linearization
\begin{align*}
L(h)&=-\frac{1}{(\partial w_0/\partial z_2)^2}
\begin{bmatrix}
	-\frac{\partial w_0}{\partial z_1} \\
	1
\end{bmatrix}
\frac{\partial}{\partial z_2}h - \frac{1}{\partial w_0/\partial z_2}\e_1\frac{\partial}{\partial z_1}h \\
&= -\frac{1}{\partial w_0/\partial z_2}\Bigl(\widetilde{\nabla} w_0 \frac{\partial}{\partial z_2}h + \e_1\frac{\partial}{\partial z_1}h \Bigr)~.
\end{align*}
If we insert this formula into the expressions for $a_1$ and $a_2$, we get that the leading order operator of the linearization of the left-hand side of \eqref{weq2} is
\begin{equation} \label{lineq}
	\sum_{i,j} \bigl[A_0\bigr]_{ij} \frac{\partial^2 h}{\partial z_i\partial z_j}
\end{equation}
with $A_0$ given by
\begin{equation} \label{Aform}
	A_0 = -\frac{1}{\partial w_0/\partial z_2} 
	\begin{bmatrix}
		\e_1 \cdot \Phi^*[I]_{H_0} \e_1 & \e_1 \cdot \Phi^*[I]_{H_0} \widetilde{\nabla} w_0 \\[1mm]
		\e_1 \cdot \Phi^*[I]_{H_0} \widetilde{\nabla} w_0 & \widetilde{\nabla} w_0 \cdot \Phi^*[I]_{H_0} \widetilde{\nabla} w_0
	\end{bmatrix}.
\end{equation}
Here $H_0^{-1}$ is defined as $H^{-1}$ with $w$ replaced by $w_0$ and $\Phi^*[I]_{H_0} = \Phi^*[I] \circ H_0^{-1}$. We similarly obtain that the highest order terms in the linearization of the operator on the left-hand side of the boundary condition \eqref{wbc2} are given by
\begin{align*} 
	&-\frac{1}{\partial w_0/\partial z_2}  \widetilde{\nabla} w_0 \cdot\bigl(2\Phi^*[I]_{H_0} \widetilde{\nabla} w_0 + (\Phi^*[I]_{H_0}-I)\e_2\bigr) \frac{\partial}{\partial z_2}h \notag \\
 	&\hspace{4cm} -\frac{1}{\partial w_0/\partial z_2}\bigl((2\widetilde{\nabla} w_0 + \e_2)\cdot\Phi^*[I]_{H_0} \e_1\bigr)\frac{\partial}{\partial z_1}h~.
\end{align*}
Since we linearize around a function $w_0$ that satisfies the boundary condition \eqref{wbc}, then this linearized boundary operator simplifies to 
\begin{equation} \label{linbc2}
	-\frac{1}{\partial w_0/\partial z_2}\widetilde{\nabla} w_0 \cdot\Phi^*[I]_{H_0} \widetilde{\nabla} w_0 \frac{\partial}{\partial z_2}h -\frac{1}{\partial w_0/\partial z_2}\bigl( (2\widetilde{\nabla} w_0+\e_2)\cdot \Phi^*[I]_{H_{0}} \e_1 \bigr)\frac{\partial}{\partial z_1}h~.
\end{equation}
A simple calculation gives that for any $\xi = (\xi_1,\xi_2)\in \R^2$,
\begin{align*}
	\xi \cdot A_0 \xi &= -\frac{1}{\partial w_0/\partial z_2}(\xi_1\e_1 + \xi_2\widetilde{\nabla} w_0) \cdot \Phi^*[I]_{H_0}(\xi_1\e_1 + \xi_2\widetilde{\nabla} w_0) \\
	&= -\frac{1}{(\partial w_0/\partial z_2)^3} \begin{bmatrix}
		\xi^\perp \cdot \nabla w_0 \\
		\xi_2
	\end{bmatrix} 
	\cdot \Phi^*[I]_{H_0} \begin{bmatrix}
		\xi^\perp \cdot \nabla w_0 \\
		\xi_2 
	\end{bmatrix}.
\end{align*}
As we assume that $w_0$ satisfies $\partial w_0/\partial z_2>0$ at $(0,0)$ then
\begin{equation} \label{lowb}
	\min_{\abs{\xi}=1} \bigl(\xi^\perp \cdot \nabla w_0(0,0)\bigr)^2 + (\xi_2)^2 > 2c \text{ for some constant } c > 0~.
\end{equation}
To see this we simply note that the function 
\begin{equation*}
	\xi \mapsto \bigl(\xi^\perp \cdot \nabla w_0(0,0)\bigr)^2 + (\xi_2)^2 
\end{equation*}
is continuous, non-negative, and only vanishes for $\xi=0$. By continuity it follows from \eqref{lowb} that
\begin{equation*}
	\min_{\abs{\xi}=1} \bigl(\xi^\perp \cdot \nabla w_0(z)\bigr)^2 + (\xi_2)^2 > c
\end{equation*}
for $z$ in a neighborhood of $(0,0)$. Due to the positive definiteness of $\Phi^*[I]_{H_0}$ it now follows immediately that 
\begin{equation} \label{ellipcond}
	-\xi \cdot A_0 \xi \geq c \abs{\xi}^2 \text{ for some constant } c > 0 
\end{equation}
for $z$ in some neighborhood of $(0,0)$. Since $\partial w_0/\partial z_2 > 0$ in a neighborhood of $(0,0)$ and $\Phi^*[I]_{H_0}$ is positive definite, it also follows that 
\begin{equation} \label{ellipbc}
 	\widetilde{\nabla} w_0 \cdot\Phi^*[I]_H \widetilde{\nabla} w_0 > c \text{ for some constant } c > 0~, 
\end{equation} 
for $z$ in some neighborhood of $(0,0)$. The bounds~\eqref{ellipcond} and~\eqref{ellipbc} ensure that the linear operator~\eqref{lineq} in combination with the oblique derivative condition~\eqref{linbc2} is uniformly elliptic and well-posed in the sense of Agmon--Douglis--Nirenberg. From~\cite[Theorem~11.2]{ADN} (with~$l=m=m_1=1$ and~$p=k+1$) we may thus infer local regularity of solutions to the non-linear problem~\eqref{weq2}--\eqref{wbc2}. To be quite precise this shows that: 
\begin{itemize}
	\item If $a_p$, $p=0,1,2$, and $b$ are in $C^{k,\alpha}$ in all arguments, $k\geq 1$, then $w$ is in $C^{k+1,\alpha}$. 
\end{itemize}
On the other hand, we already know that:
\begin{itemize}
	\item If $\Phi^*[I]$ is in $C^{k+1,\alpha}$ for some $k\geq 1$, then $\Phi^*[I]_H$ is in $C^{k+1,\alpha}$ and $\frac{\partial}{\partial z_2} \Phi^*[I]_H$ is in $C^{k,\alpha}$ in all arguments\footnote{Notice that the arguments of $\Phi^*[I]_H$ and $\frac{\partial}{\partial z_2} \Phi^*[I]_H$ are $(z_1,w)$,  and $(z_1,w,\frac{\partial w}{\partial z_2})$, respectively}. Consequently $a_p$, $p=0,1,2$, and $b$ are in $C^{k,\alpha}$ in all arguments.
\end{itemize}
A combination of these two assertions shows that $w$ and therefore the boundary $\partial\Omega$ is of class $C^{k+1,\alpha}$ near $P = (0,0)$. This completes the proof in the remaining case of Theorem~\ref{main}.

\begin{remark}
	The same approach as above would also work with Assumption~\ref{assump:NonDeg4} for any vector $\xi$ (not just for the normal vector) --- in that case one would just work with the function $v(y)=\xi^\T \Phi^{-1}(y)-\xi^\T y$ in place of $(\Phi^{-1}(y))_n-y_n$. However, as we already saw in the previous section, when $\xi$ has a tangential component the added complexity of introducing the Hodograph transform is not necessary.
\end{remark}

\subsection*{Acknowledgements}

HG was partially supported by grant 10.46540/3120-00003B from Independent Research Fund Denmark. MSV was partially supported by NSF grant DMS-22-05912.

\bibliographystyle{plain}

\end{document}